\documentclass[11pt]{article}
\usepackage{amssymb, amsfonts, amsmath, amsthm}
\usepackage{color}
\usepackage{hyperref, cite}
\usepackage{authblk}

\newtheorem{theorem}{Theorem}[section]

\newtheorem{prop}[theorem]{Proposition}

\numberwithin{equation}{section}

\def\pf{\noindent {\it Proof.\ }}

\begin{document}

\title{\textbf{A construction method for WZ seeds}}
\author{Qing-Hu Hou$^{1,}$\footnote{Supported by the National Key Research and Development Program of China (2023YFA1009401).} \ and Yan-Ping Mu$^{2,}$\footnote{Corresponding author}
	\footnote{Supported by the National Natural Science Foundation of China (grant 12171255).}}
\affil{$^1$School of Mathematics and KL-AAGDM, Tianjin University,\\ Tianjin 300350, China \\
\texttt{qh\_hou@tju.edu.cn} \\ \vspace{5pt}
$^{2}$School of Mathematical Sciences, Tianjin University of Technology,\\ Tianjin 300384, China \break
\texttt{yanping.mu@gmail.com}}
\maketitle

\begin{abstract}
We propose a systematic method for constructing Wilf-Zeilberger (WZ) seeds and present seven WZ seeds. We also demonstrate how to construct WZ seeds from existing ones. With these WZ seeds,  several hypergeometric identities are derived. The construction can be extended to the $q$-cases, leading to the $q$-analogues of the seven WZ seeds.

\end{abstract}



\section{Introduction}
Recall that a function $F(n_1, \ldots,n_r)$ is  {\it hypergeometric} in the variables $n_1,\ldots,n_r$ if all the ratios 
\[
\frac{F(n_1+1,n_2,\ldots,n_r)}{F(n_1,n_2,\ldots,n_r)},\ \frac{F(n_1,n_2+1,n_3,\ldots,n_r)}{F(n_1,n_2,\ldots,n_r)}, \ \ldots, \frac{F(n_1,\ldots,n_{r-1}, n_r+1)}{F(n_1,n_2,\ldots,n_r)}
\]
are rational functions of $n_1,n_2,\ldots,n_r$.  A pair of hypergeometric terms $F(n,k)$ and $G(n,k)$ is called a {\it Wilf-Zeilberger {\rm (}WZ\/{\rm )} pair}  \cite{WZ92, Pet96} if  
\begin{equation}\label{WZ}
\Delta_n F(n,k) = \Delta_k G(n,k),
\end{equation}
where $\Delta_n f(n) = f(n+1) - f(n)$ is the difference operator with respect to $n$. In this case, we refer to $G(n,k)$ as the {\it WZ mate} of $F(n,k)$.

Gessel \cite{Gess94}  provided a systematic method to generate WZ pairs. Specifically, given a hypergeometric summation formula with parameters $a,b,\ldots$
\[
\sum_k F(k, a,b,\ldots) = H(a,b,\ldots),
\]
the term 
\[
F_0(k,a,b,\ldots) = \frac{F(k,a,b,\ldots)} {H(a,b,\ldots)}
\]
often possesses the following property (we call it the {\it universal property}): For any integers $K,A,B,\ldots$ and any parameter $k_0$, 
\begin{equation*}
F(n,k) =F_0 (Kn+k_0+k, An+a, Bn+b,\ldots)
\end{equation*}
has a WZ mate. By selecting suitable $K,A,B,\ldots$ and $k_0,a,b,\ldots$, Gessel derived a large number of terminating hypergeometric identities.

Guillera \cite{Gu06, Gu10} provided a method for  constructing WZ pairs with given forms of $G(n,k)$. Ding-Mu \cite{DM} generalized the forms of $G(n,k)$ and constructed WZ pair by the extended Zeilberger algorithm.

A WZ pair can be also obtained from a hypergeometric term $F(n,k)$ satisfying
the telescoped recurrence
\begin{equation}\label{tele}
p_0(n) F(n,k) + p_1(n) F(n+1,k) = \Delta_k G(n,k).
\end{equation}
In fact, 
\[
F(n,k) \cdot  (-1)^n \prod_{j=0}^{n-1} \frac{p_1(j)}{p_0(j)}
\]
has a WZ mate. Following this approach, Campbell \cite{Camp25}  derived $q$-analogues of $\pi$-formulas due to Ramanujan and Guillera, and Ding-Mu \cite{DM} obtained $q$-analogues of the identity (H.1) due to Van Hamme.

Recently, Au \cite{Au25a, Au25b} named a hypergeometric term $F_0(k,a,b,\ldots)$ satisfying the universal property a {\it WZ seed} in variables $a,b,\ldots$ and listed $13$ WZ seeds, most of which are due to Gessel. By selecting suitable WZ seeds and parameters, he presented many hypergeometric identities and confirmed several conjectures posed by Sun \cite{Sun24, Sun26}. 

Inspired by the telescoped recurrence \eqref{tele},  in Section 2 we present a construction method for WZ seeds starting with appreciate hypergeometric terms $F(k,a,b,\ldots)$. If for each parameter $a,b,\ldots$, the Zeilberger algorithm \cite{Zeil90, Zeil91} yields a telescoped recurrence involving only two terms, then we can modify $F(k,a,b,\ldots)$ to obtain a WZ seed. We note that all WZ seeds introduced by Au can be recovered by this method. We further provide seven WZ seeds, some of which generalize those introduced by Au with one more parameter.  Additionally, Zeilberger \cite{Pet96} and Gessel \cite{Gess94} have observed that if $(F(n,k), G(n,k))$ is a WZ pair, then so is $(G(k, n), F(k, n))$. We show that WZ seeds exhibit a similar phenomenon, which enables us to construct WZ seeds from existing ones.

As applications, we establish several hypergeometric identities in Section~3. We first use the second WZ seed to illustrate evaluating hypergeometric series. Then we utilize the first WZ seed to derive some transformation formulas. At last, we present one or two examples for each of the remaining WZ seeds.

In Section 4, we extend the construction to the $q$-cases and give the $q$-analogues of the seven WZ seeds, which also generalize those given by Au \cite{Auc}.  We also discuss the change from the base $q$ to $q^{-1}$, which leads to new $q$-WZ seeds.

\section{Constructing  WZ seeds}
Let $F(k, a,b,\ldots)$ be a hypergeometric term  with parameters $a,b,\ldots$.
Suppose that when we apply the Zeilberger algorithm to $F$ with shifts in $a$ and obtain the telescoped recurrence involving only two terms
\[
p_0(a,b,\ldots) F(k, a,b,\ldots) + p_r(a,b,\ldots) F(k, a+r,b,\ldots) = \Delta_k G(k,a,b,\ldots),
\]
where $p_0$ and $p_r$ are polynomials in $a$ and independent of $k$. Factoring $p_0$ and $p_r$, we can write
\begin{equation}\label{p0pr}
p_0 = c_0 (a/r-u_1)  \ldots (a/r-u_l), \quad p_r = c_r (a/r-v_1)  \ldots (a/r-v_m).
\end{equation}
Then we set
\begin{equation}\label{F2Fa}
F_a(k,a,b,\ldots) =   \left(-\frac{c_r}{c_0} \right)^a \frac{\prod_{j=1}^m \Gamma(a-v_j)}{\prod_{i=1}^l \Gamma(a-u_i)} \cdot F(k, r \cdot a, b,\ldots) 
\end{equation}
and
\begin{equation}\label{Ga}
G_a(k,a,b,\ldots) =  \frac{c_r^a}{(-c_0)^{a+1}} \frac{\prod_{j=1}^m \Gamma(a-v_j)}{\prod_{i=1}^l \Gamma(a-u_i+1)} \cdot G(k,a \cdot r,b,\ldots),
\end{equation}
where $\Gamma(x)$ denotes the Gamma function of $x$.
It is routine to check that 
\begin{equation}\label{Fa}
\Delta_a F_a(k, a,b,\ldots) = \Delta_k G_a(k,a,b,\ldots).
\end{equation}

Now we apply the Zeilberger algorithm to $F_a(k,a,b,\ldots)$ with shifts in $b$. Suppose that 
we obtain a telescoped recurrence of the form
\begin{multline}\label{Fap}
p'_0(a,b,\ldots) F_a(k, a,b,\ldots) + p'_{r'}(a,b,\ldots) F_a(k, a, b+r',\ldots) \\
= \Delta_k G'(k,a,b,\ldots),
\end{multline}
where $p'_0$ and $p'_{r'}$ are polynomials in $b$ and independent of $k$. We remark that 
in most cases, $p'_0/p'_{r'}$ is also independent of $a$. One explanation may be as follows. Suppose
\[
\sum_{k = - \infty}^{\infty} F_a(k,a,b,\ldots)
\] 
converges to $f(a,b,\ldots)$
and 
\[
\lim_{k \to \pm \infty} G(k,a,b,\ldots) = \lim_{k \to  \pm \infty} G'(k,a,b,\ldots) = 0. 
\]
Summing both sides of \eqref{Fa} and \eqref{Fap} over $k$ from $-\infty$ to $+\infty$, we get
\[
f(a,b,\ldots) = f(a+1,b,\ldots)
\]
and
\[
p'_0(a,b,\ldots) f(a,b,\ldots)  = -p'_{r'}(a,b,\ldots) f(a, b+r',\ldots),
\]
implying that 
\[
\frac{p'_0(a,b,\ldots)}{p'_{r'}(a,b,\ldots)} = \frac{p'_0(a+1,b,\ldots)}{p'_{r'}(a+1,b,\ldots)}.
\]

Given the telescoped recurrence  \eqref{Fap}, we can factor $p'_0$ and $p'_{r'}$ and construct $F_b$ similar to \eqref{F2Fa} such that $\Delta_b F_b$ has a WZ mate. Since $p'_0/p'_{r'}$ is independent of $a$, $\Delta_a F_b$ still has a WZ mate. As shown by Au \cite[Theorem 3.3]{Au25a},  $F_0(k,a,b,\ldots)$ is a WZ seed if and only if there exist hypergeometric terms  $G_a, G_b, \ldots$ such that
\[
\Delta_a  F_0(k,a,b,\ldots) = \Delta_k G_a, \quad  \Delta_b F_0(k,a,b,\ldots)= \Delta_k G_b, \quad  \ldots.
\]
By repeating this process to other parameters, we finally obtain a WZ seed $F_0(k, a,b, \ldots)$.

Let us illustrate the above process with 
\[
F(k,a,b,c,d) = \frac{(a)_k (b)_k}{(c)_k (d)_k},
\]
where $(a)_k = \Gamma(a+k)/\Gamma(a)$ is the rising factorial  (also known as the Pochhammer symbol).
For simplify, we will abbreviate the variables that do not change.
By the Zeilberger algorithm, we find that
\begin{multline*}
(a-d+1) (a-c+1) F(a)  -(a+b-c-d+2) a F(a+1) \\
= \Delta_k \big( -(c+k-1)(d+k-1) F(a) \big).
\end{multline*}
Now set
\[
F_a =  \Gamma\! \left[ \begin{array}{c} a,\ a+b-c-d+2 \\ a+1-c,\ a+1-d \end{array} \right] \cdot F.
\]
Here and in what follows, we use the shorthand notation
\[
\Gamma\! \left[ \begin{array}{c}a_1, \ldots, a_r \\ b_1, \ldots, b_s \end{array} \right] = \frac{\prod_{i=1}^r \Gamma(a_i)}{\prod_{j=1}^s \Gamma(b_j)}.
\]

Also by the Zeilberger algorithm, we find that 
\begin{multline*}
(b-d+1) (b-c+1) F_a(b) - b F_a(b+1) \\
=  \Delta_k \big( -(c+k-1)(d+k-1) F_a(b) \big).
\end{multline*}
Hence we set
\[
F_b =  \Gamma\! \left[ \begin{array}{c} b \\ b+1-c,\ b+1-d \end{array} \right]   \cdot F_a.
\]
Applying the Zeilberger algorithm again, we get
\[
c F_b(c) + F_b(c+1) = \Delta_k \left( \frac {c \left( d+k-1 \right) }{-c-d+a+b+1} F_b(c) \right),
\]
so that
\[
F_c =  \frac{(-1)^c}{\Gamma(c)} {F_b}.
\]
Finally, we have
\[
d F_c(d) + F_c(d+1) = \Delta_k \left( \frac {d \left( c+k-1 \right) }{-c-d+a+b+1} F_c(d) \right).
\]
This leads to the second WZ seed in the following theorem. With the same method, we construct seven WZ seeds.

\begin{theorem}\label{th-WZ}
The following functions are WZ seeds:
\begin{itemize}
\item[{\rm (1)}]
\begin{equation}\label{WZ-abz}
F_0(k,a,b) = {z}^{b+k} \left( z-1 \right) ^{a-b} \ \Gamma\! \left[ \begin{array}{c} 
	a + k \\ 
	b + k 
\end{array} \right] \Gamma  \left( -a+b	\right),
\end{equation}

\item[{\rm (2)}] 
\begin{multline}\label{Sabcd}
F_0(k,a,b,c,d) = (-1)^{c+d} \ \Gamma\! \left[ \begin{array}{c} a+k,\ b+k \\ c+k,\ d+k \end{array} \right] \\
\times \Gamma\! \left[ \begin{array}{c}  a+b-c-d+2 \\ a-c+1,\ a-d+1,\ b-c+1,\ b-d+1 \end{array} \right] .
\end{multline}

\item[{\rm (3)}]
\begin{multline}\label{WZ-abc'}
F_0(k,a,b,c) =  \left( -1 \right) ^{k} \ \Gamma\! \left[ \begin{array}{c} 
2a + k,\ 2b + k \\ 
2a + 2c + k,\ 2b + 2c + k 
\end{array} \right] \\[5pt]
\times \Gamma\! \left[ \begin{array}{c} 
2c,\ a  - b + c +\frac{1}{2} \\ 
a - b - c+\frac{1}{2} 
\end{array} \right],
\end{multline}

\item[{\rm (4)}]
{\small
\begin{multline}\label{WZ-abc11}
F_0(k,a,b,c,d,e) =  \left( -1 \right) ^{
d+e} \ \Gamma\! \left[ \begin{array}{c} 
a+k,\ b+k,\ c+k \\ 
d+k,\ e+k,\ a + b + c - d - e + 2 + k 
\end{array} \right]
\\[5pt]
\times	\ \Gamma\! \left[ \begin{array}{c} 
a + b - d - e + 2,\ a + c - d - e + 2,\ b + c - d - e + 2 \\ 
a - d + 1,\ a - e + 1,\ b - d + 1,\ b - e + 1,\ c - d + 1,\ c - e + 1 
\end{array} \right],
\end{multline}}

\item[{\rm (5)}]
$F_0(k,a,b,c,d) =$
{\small
\begin{multline}\label{WZS4}
\ \Gamma\! \left[ \begin{array}{c} 
	a + k,\ b + k,\ c + k \\ 
	2d - a + k,\ 2d - b + k,\ 2d - c + k 
\end{array} \right]
\\[5pt]
\times \ \Gamma\! \left[ \begin{array}{c} \begin{array}{c}
a + b + c - 3d+\frac{3}{2} \\
a + b - 2d + 1,\ a +c  - 2d  + 1,\ b + c - 2d + 1 \end{array}
\end{array} \right]
\\[5pt]
\times \frac{ \left( -1 \right) ^{a+b+c+d}}{\ \Gamma\! \left[ \begin{array}{c} 
a - d+\frac{1}{2},\ b - d+\frac{1}{2},\ c - d+\frac{1}{2} 
\end{array} \right] },
\end{multline}}

\item[{\rm (6)}]
$F_0(k,a,b,c,d,e) =  \left( -1 \right) ^{a+b+c+d+e}  \left( e+2\,k \right)$
{\small 
\begin{multline}\label{WZS5}
\times 
\ \Gamma\! \left[ \begin{array}{c} 
a + k,\ b + k,\ c + k,\ d + k \\ 
e - a + 1 +  k ,\  e - b  + 1 +  k,\  e - c  + 1 +  k,\ e - d  + 1 +  k 
\end{array} \right]
\\[5pt]
\times \ \Gamma\! \left[ \begin{array}{c} 
a + b + c + d - 2e \\ 
a + b - e,\ a + c  - e,\ a + d - e,\ b + c - e,\ b + d - e,\ c + d - e 
\end{array} \right],
\end{multline}}

\item[{\rm (7)}]
$F_0(k,a,b,c,d,e,f) = \left( -1 \right) ^{f} \left( f+2\,k \right) $
{\small 
\begin{multline}\label{WZS6}
\times	\ \Gamma\! \left[ \begin{array}{c} 
a + k,\ b + k,\ c + k\\ 
f - a  + 1 + k,\   f - b  + 1 + k,\  f - c  + 1 + k 
\end{array} \right]
\\[5pt]
\times	\ \Gamma\! \left[ \begin{array}{c} 
	d + k,\ e + k,\ 3 f-a-b-c-d-e+1+k \\ 
 f - d  + 1 + k,\ f - e  + 1 + k,\  a+b+c+d+e-2f+k 
\end{array} \right]
\\[5pt]
	\times\Gamma \left[ \begin{array}{c} 
		a + b + c + d - 2f,\ a + b + c + e - 2f,\ a + b + d + e - 2f \\ 
		a + b - f,\ a + c - f,\ a + d - f,\ a + e - f,\ b + c - f,\ b + d - f 
	\end{array} \right]
 \\[5pt]
	\times \Gamma \left[ \begin{array}{c} 
		a + c + d + e - 2f,\ b + c + d + e - 2f  \\ 
		b + e - f,\ c + d - f,\ c + e - f,\ d + e - f
	\end{array} \right].
\end{multline}}

\end{itemize}
\end{theorem}
\pf
We only list the hypergeometric terms which we start with.
\begin{itemize}
\item[{\rm (1)}]
$\frac{(a)_k}{(b)_k} z^k$,

\item[{\rm (2)}]
$\frac{(a)_k (b)_k}{(c)_k (d)_k}$,

\item[{\rm (3)}]
$\frac{(a)_k (b)_k}{(c+a)_k (c+b)_k} (-1)^k$,

\item[{\rm (4)}]
$\frac{(a)_k (b)_k (c)_k}{(d)_k (e)_k (2+a+b+c-d-e)_k}$,

\item[{\rm (5)}]
$\frac{(a)_k (b)_k (c)_k}{(d-a)_k (d-b)_k (d-c)_k}$,

\item[{\rm (6)}]
$(e+2k) \frac{(a)_k (b)_k (c)_k (d)_k}{(1+e-a)_k (1+e-b)_k (1+e-c)_k (1+e-d)_k}$,

\item[{\rm (7)}]
$(f+2k) \frac {(a)_k (b)_k (c)_k (d)_k (e)_k (1+3f-a-b-c-d-e)_k}{(1+f-a)_k (1+f-b)_k (1+f-c)_k (1+f-d)_k(1+f-e)_k (a+b+c+d+e-2f)_k}$.
\qed
\end{itemize}

\noindent {\bf Remark.} The WZ seed \eqref{Sabcd} can be derived from the Dougall's summation formula \cite[Theorem2.8.2]{AAR}: for $1 + \operatorname{Re}(a + b) < \operatorname{Re}(c + d)$,
\[
\sum_{k=-\infty}^{\infty} \frac{\Gamma(a + k)\Gamma(b + k)}{\Gamma(c + k)\Gamma(d + k)} 
= \frac{\pi^2}{\sin \pi a \, \sin \pi b} \ \Gamma\! \left[ \begin{array}{c} 
c + d - a - b - 1 \\ 
c - a,\ d - a,\ c - b,\ d - b 
\end{array} \right],
\]
with the help of Euler’s reflection formula
\[
\Gamma(x) \Gamma(1-x) = \frac{\pi}{\sin (\pi x)}.
\]

The WZ seeds in Theorem~\ref{th-WZ} generalize those given in \cite{Au25a}. Specifically:
\begin{itemize}
\item 
By setting $d=1$, \eqref{Sabcd} is equivalent to {\tt Gauss2F1}.
\item 
By setting $e=a+b+c+1-d$, \eqref{WZ-abc11}  is equivalent to {\tt Balanced3F2}.
\item 
By setting $d=a+1$, \eqref{WZS4} is equivalent to {\tt Dixon3F2}.
\item 
By setting $e=a$, \eqref{WZS5}  is equivalent to {\tt Dougall5F4}.
\item 
By setting $f=a$, \eqref{WZS6} is equivalent to {\tt Dougall7F6}.
\end{itemize}


Zeilberger \cite{Pet96} and Gessel \cite{Gess94} have observed that if $(F(n,k), G(n,k))$ is a WZ pair, then so is $(G(k, n), F(k, n))$. We find that WZ seeds exhibit a similar phenomenon.
\begin{theorem}
Let $F_0(k,a,b,\ldots)$ be a WZ seed such that for any two parameters, say $a$ and $b$, $\Delta_a \Delta_b F_0$ is not a rational function of $k$. Let
\[
F(n,k) = F_0(Kn+k, An+a, Bn+b,\ldots),
\]
where $K,A,B,\ldots$ are integers, and let $G(n,k)$ be the WZ mate of $F(n,k)$. Suppose $G(n,k)$ is not a rational function of $k$. Then
$G(k,n)$ is a WZ seed with parameters $n, a,b,\ldots$.
\end{theorem}

\pf It is known that
\[
\Delta_n G(k,n) = \Delta_k F(k,n).
\]
By the symmetry of $a,b,\ldots$, it suffices to show that  there exists a hypergeometric term $H(k)$ such that
\[
\Delta_a G(k,n) = \Delta_k H(k).
\]

We first give an explicit expression for $G$. Since $F_0(k,a,b,\ldots)$ is a WZ seed, we can assume
{
\small
\[
\Delta_a F_0(k,a,b,\ldots)  = \Delta_k G_a(k,a,b,\ldots) , \ \Delta_b F_0(k,a,b,\ldots)  = \Delta_k G_b(k,a,b,\ldots) , \ \ldots.
\]
}
It is routine to verify that
\begin{align*}
\Delta_n F(n,k)  = & F_0(Kn+K+k,An+A+a,Bn+B+b,\ldots) \\
& \quad - F_0(Kn+k,An+a,Bn+b,\ldots)  \\
= & \sum_{i=0}^{K-1} \Delta_k F_0(Kn+k+i, An+A+a,Bn+B+b, \ldots) \\
& + \sum_{i=0}^{A-1} \Delta_a F_0(Kn+k, An+a+i,Bn+B+b, \ldots) \\
& + \sum_{i=0}^{B-1} \Delta_b F_0(Kn+k, An+a,Bn+b+i, \ldots) + \cdots \\
=& \Delta_k \left( \sum_{i=0}^{K-1} F_0(Kn+k+i, An+A+a,Bn+B+b, \ldots) \right. \\
& + \sum_{i=0}^{A-1} G_a (Kn+k, An+a+i,Bn+B+b, \ldots) \\
& \left. + \sum_{i=0}^{B-1} G_b (Kn+k, An+a,Bn+b+i, \ldots) + \cdots \right) \\
:= & \Delta_k \tilde{G}(n,k).
\end{align*}
We thus derive that 
\[
G(n,k) -  \tilde{G}(n,k)  = C
\]
is a constant  with respect to $k$. 
Recall that, according to \cite[Proposition 5.6.2]{Pet96},  the sum of two hypergeometric terms is a non-zero hypergeometric term if and only if they are similar. That is,  the ratio of these two terms is a rational function of $k$. If $C \not= 0$, then $G$ is similar to $C$. This implies that $G$ is a rational function of $k$, which contradicts to the hypothesis. Therefore, we conclude that
$G=\tilde{G}$.

Next, notice that 
\[
\Delta_k  \Delta_a G_b = \Delta_a \Delta_k G_b = \Delta_a \Delta_b F_0 =  \Delta_b \Delta_a F_0 = \Delta_k  \Delta_b G_a,
\]
i.e., 
\[
\Delta_k (\Delta_a G_b - \Delta_b G_a) = 0.
\]
Since $\Delta_k  \Delta_a G_b = \Delta_a \Delta_b F_0$ is not a rational function of $k$,  neither is $  \Delta_a G_b$.  
Similar to the above discussion, we deduce that $\Delta_a G_b = \Delta_b G_a$.

Now,
\begin{align*}
\Delta_a G(k,n)  = & \Delta_a  \left( \sum_{i=0}^{K-1} F_0(Kk+n+i, Ak+A+a,Bk+B+b, \ldots) \right. \\
& \qquad + \sum_{i=0}^{A-1} G_a (Kk+n, Ak+a+i,Bk+B+b, \ldots) \\
& \qquad \left. + \sum_{i=0}^{B-1} G_b (Kk+n, Ak+a,Bk+b+i, \ldots) + \cdots \right) \\
= & \sum_{i=0}^{K-1} \big( G_a(Kk+n+i+1, Ak+A+a,Bk+B+b, \ldots) \\
& \qquad  - G_a(Kk+n+i, Ak+A+a,Bk+B+b, \ldots) \big) \\
& + \sum_{i=0}^{A-1}  \Delta_a G_a (Kk+n, Ak+a+i,Bk+B+b, \ldots) \\
& + \sum_{i=0}^{B-1}  \Delta_b G_a (Kk+n, Ak+a,Bk+b+i, \ldots) + \cdots \\
=  & \Delta_k G_a(Kk+n, Ak+a,Bk+b, \ldots). \tag*{\qed}
\end{align*}

Let us illustrate the above construction with a simple example.  By setting $a \to n+a$ and $b \to 2n+b$ in  \eqref{WZ-abz}, we get
\[
F(n,k) = {z}^{b+2n+k} \left( z-1 \right) ^{a-b-n} \ \Gamma\! \left[ \begin{array}{c} 
a + n + k  \\ 
b + 2n +k 
\end{array} \right] \Gamma  \left( -a+b+n	\right).
\]
This leads to the WZ seed
\[
((2-z)k+az-bz+b+n)
z^{b+n+2k}  (z-1)^{a-b-k-1} 
\ \Gamma\! \left[ \begin{array}{c} 
a + n + k, \ -a + b +k \\ 
b+2k+n+1
\end{array} \right] ,
\]
which is equivalent to 
\[
((2-z)k+az-bz+b+z-1) z^{b-1+2k} (z-1)^{a-b-k}  \ \Gamma\! \left[ \begin{array}{c} 
a  + k, b-a+k-1  \\ 
b  + 2 k 
\end{array} \right].
\]

\section{Deriving identities}
Given a WZ seed $F_0(k,a,b,\ldots)$, one can choose any integers $K,A,B,\ldots$ and arbitrary parameters $k_0,a,b,\ldots$ to obtain
\[
F(n,k) = F_0(Kn+k_0+k, An+a, Bn+b,\ldots)
\]
and its WZ mate $G(n,k)$. By  Proposition 2.1 in \cite{Au25a}, we have
\begin{equation}\label{FG}
	\sum_{n=0}^\infty G(n,0) = \sum_{k=0}^\infty F(0,k)  + \sum_{n=0}^\infty g(n)  - \lim_{n \to \infty} \sum_{k=0}^\infty F(n,k),
\end{equation}
where
\[
g(n) =  \lim_{k \to \infty} G(n,k).
\]
In most cases, $g(n)$ and $\lim\limits_{n \to \infty} \sum\limits_{k=0}^\infty F(n,k)$ vanish.

We begin by employing the WZ seed \eqref{Sabcd} to evaluate certain hypergeometric series. It is worth noting that one may substitute
$\Gamma(a)$ with $(-1)^a/ \Gamma(1-a)$ and the resulting expression remains a WZ seed. For instance, we shall utilize the following variant of the WZ seed \eqref{Sabcd}:
\begin{multline*}
F_0(k,a,b,c,d) =(-1)^b {\frac {\Gamma  \left( -b+d \right) \Gamma  \left( -b+c \right) 
		\Gamma  \left( -a+d \right) \Gamma  \left( -a+c \right) }{\Gamma 
		\left( d+k \right) \Gamma  \left( c+k \right) }} \\
\times
\frac{\Gamma 
	\left( b+k \right) \Gamma  \left( a+k \right) \Gamma  \left( a
	\right) \Gamma  \left( 1-a \right) } {\Gamma  \left( -1-a-b+c
	+d \right) }.
\end{multline*}

To obtain the evaluation, we select parameters $k_0,a,b,c,d$ such that $\sum_{k=0}^\infty F(0,k)$ has a closed form expression. A suitable choice is given by  
\[
k_0=0,a=\frac{1}{3},b=1,c= \frac{4}{3}, d=2, 
\]
resulting in  
\[
F(0,k) = - \frac {8{\pi }^{2}}{ 3 \left( 1+k \right)  \left( 1+3\,k \right) }.
\]
By {\tt Maple} or {\tt WolframAlpha}, we find that
\begin{equation} \label{S0}
S_0= \sum_{k=0}^\infty \frac{1}{(3k+1)(k+1)} =   \frac{\pi \sqrt{3}}{12} + \frac{3}{4} \log 3.
\end{equation}
By selecting distinct parameters $K,A,B,C,D$, we generate a family of infinite series whose evaluations are all expressed in terms of $S_0$.
\begin{prop}
Let $S_0$ be given by \eqref{S0}.
Then
\begin{align}
& \label{Ex1-cd}	\sum_{n=0}^\infty  {\frac {  (9\,n+7)   }{ \left( n+1 \right) 
		\left( 3\,n+1 \right)   }} \frac{\left(\frac{5}{3}\right)_n}{\left(\frac{3}{2}\right)_n}  \frac{1}{4^n}
= 6 S_0. \\
& \label{Ex1-kc}
\sum_{n=0}^\infty  \frac{ 9\,n+5 }{  n+1 }  \frac{\left(\frac{1}{3}\right)^2_n}{\left(\frac{5}{3}\right)_n \left(\frac{7}{6}\right)_n}  \frac{1}{4^n}
= 4S_0. \\
&\label{Ex1-kcd}
\sum_{n=0}^\infty \left( 189\,{n}^{2}+222\,n+61 \right)  \frac {\left(\frac{1}{3}\right)^2_n (1)_n}{ \left(\frac{7}{6}\right)_n \left(\frac{3}{2}\right)_n^2}   \frac {1}{64^n} = 48S_0. \\
& \label{Ex1-kbcd}
\sum_{n=0}^\infty \frac{30 n^2+35 n+11}{(1+2 n) (1+n)} \cdot  \frac{\left(\frac{1}{3}\right)_n}{\left(\frac{7}{6}\right)_n} \frac{1}{(-4)^n}
=  8S_0. \\
& \label{Ex1-a2cd}
\sum_{n=0}^\infty \frac{90 n^2+111 n+31}{(1+6 n) (2+3 n) (1+n)} \cdot  \frac{\left(\frac{1}{3}\right)_n}{\left(\frac{3}{2}\right)_n} \frac{1}{(-4)^n}
=  12S_0.
\end{align}

\end{prop}
\pf 
Using the parameters $K,A,B,C,D$ given in Table \ref{t-1} and applying \eqref{FG}, we obtain the required identities. \qed
\begin{table}[ht]
\begin{center}
\begin{tabular}{c|ccccc}
Equation number & K & A & B & C & D  \\ \hline 
\eqref{Ex1-cd} & 0 & 0 & 0 & 1 & 1\rule{0pt}{15pt} \\[5pt] 
\eqref{Ex1-kc} & 1 & 0 & 0 & 1 &0 \\[5pt] 
\eqref{Ex1-kcd} & 1 & 0 & 0 & 1 &1 \\[5pt] 
\eqref{Ex1-kbcd} & 1 & 0 & 1 & 1 &1 \\[5pt] 
\eqref{Ex1-a2cd} & 0 & 1 & 0 & 2 &1 
\end{tabular}
\end{center}
\caption{The values for the parameters $K,A,B,C,D$. \label{t-1}}
\end{table}

Setting $k_0=\frac{1}{3}$, keeping $a,b,c,d$ unchanged, and taking the same $K,A$, $B,C,D$ as in the above proof, we obtain another family of infinite series.
\begin{prop}
Let
\[
S_1 = \sum_{k=0}^\infty \frac{1}{(3k+2)(3k+4)} = \frac{1}{2} - \frac{\sqrt{3}\pi}{18}.
\]
We have
\begin{align}
	& \sum_{n=0}^\infty  \frac { 1}{ (3n+1) (3n+4) {2n+2 \choose n+1}} 
	= \frac{2}{3} S_1, \label{Ex2-cd}\\[5pt]
	& \sum_{n=0}^\infty \frac { \left( 27\,{n}^{2}+51\,n+23 \right) \left( \frac{2}{3} \right)_n}{ \left( 3\,n+1 \right)  \left( 
		3\,n+4 \right)  \left( n+1 \right)   \left( \frac{11}{6} \right)_n} \frac{1}{4^n} = 30 S_1, \label{Ex2-kc} \\[5pt]
	& \sum_{n=0}^\infty \frac{189 n^3+495 n^2+414 n+110}{(6 n+1) (6 n+7) (6 n+5) (2 n+1)} \frac{1}{{6n \choose 3n}} = 16 S_1, \label{Ex2-kcd} \\[5pt]
	& \sum_{n=0}^\infty \frac{270n^3+720n^2+633n+185}{(n+1)(6n+7)} \frac{(\frac{2}{3})^2_n}{(\frac{4}{3})_n (\frac{11}{6})_n} \frac{1}{(-4)^n} = 120 S_1. \\[5pt]
	& \sum_{n=0}^\infty
	(30 n^2+45 n+16) \frac{(\frac{1}{6})_n (\frac{2}{3})^2_n (1)_n}{(\frac{4}{3})_n (\frac{7}{3})_n (\frac{3}{2})_n (\frac{11}{6})_n} \frac{1}{(-4)^n} 
	= 80 S_1.
\end{align}
\end{prop}

\noindent {\bf Remark.}
By the extended Zeilberger algorithm \cite{CHM}, we find that
\[
\frac{1}{{2n \choose n}} + 2 \frac { 1}{ (3n+1) (3n+4) {2n+2 \choose n+1}}  = \Delta_n \left( 
- \frac{2 (2n+1) }{(3n+1) {2n \choose n}}
\right).
\]
Therefore, \eqref{Ex2-cd} is equivalent to the known formula
\begin{equation}\label{1/bino(2n,n)}
\sum_{n=0}^\infty \frac{1}{{2n \choose n}} = \frac{2}{27} (18+\sqrt{3} \pi).
\end{equation}
Similarly, 
\[
\frac{1}{12} \frac{189 n^3+495 n^2+414 n+110}{(6 n+1) (6 n+7) (6 n+5) (2 n+1)} \frac{1}{{6n \choose 3n}} - \sum_{\ell=0}^2 \frac{1}{{2(3n+\ell) \choose 3n+\ell}} = \Delta_n z_n,
\]
for some hypergeometric term $z_n$. Therefore, \eqref{Ex2-kcd} is also equivalent to  \eqref{1/bino(2n,n)}.

Starting from the WZ seed \eqref{WZ-abz}, we derive the following transformation formulas.
\begin{prop}
Given that the following series are convergent, we have
\begin{align}
	\sum_{k=0}^\infty \frac{(a)_k}{(b)_k} z^k & = \frac{1}{1-z} \sum_{n=0}^\infty \frac{(b-a)_n}{(b)_n} \left( \frac{z}{z-1} \right)^n \label{abz-b}\\
	&  = \frac{1}{1-z} \sum_{n=0}^\infty \frac{((a-b-n)z+b+2\,n)}{ \left( 2\,n+b \right) }  \frac{(a)_n (b-a)_n}{(\frac{b+1}{2})_n (\frac{b}{2})_n} \left( \frac{z^2}{4(z-1)} \right)^n \label{abz-bk} \\
	& = \left( a-b \right) \sum_{n=0}^\infty {\frac { \left( (a+2n)z-a+b-n-1 \right)   }{
			\left( a+n-b \right)  \left( a+n+1-b \right) }} 
	\frac{(\frac{a+1}{2})_n (\frac{a}{2})_n}{(a-b)_n (b)_n} \left(4 z(1-z) \right)^n \label{abz-ak}\\	
	& = \frac{1}{1-z} \sum_{n=0}^\infty \frac{p(n)}{(1+3 n+b) (3 n+b)} \frac{(b-a)_n (\frac{a+1}{2})_n (\frac{a}{2})_n}{(\frac{b}{3})_n (\frac{b+1}{3})_n (\frac{b+2}{3})_n} \left( \frac{4z^3}{27(z-1)}\right)^n, \label{abz-bk2}
\end{align}
where 
{\small 
	\begin{multline*}
		p(n) = a^2 z^2-abz^2+anz^2-2bnz^2-2n^2z^2+abz+3anz-b^2z-4bnz-3n^2z\\
		+az+b^2+6bn-bz+9n^2-nz+b+3n.
	\end{multline*}
}
\end{prop}
\pf By setting $(a,b,k)$ to be $(a,b+n,k)$, $(a,b+n,k+n)$,  $(a+n,b,k+n)$ and $(a,b+n,k+2n)$ in \eqref{WZ-abz}, respectively, and using \eqref{FG}, we derive the desired identities. \qed

Similarly, by first replacing $b$ and $k$ with $1+n$ and $1+k+Kn\, (1 \le K \le 3)$ respectively,  and then taking the limit as $a \to 0$, we can derive the following summation formulas.
\begin{prop}
We have
\begin{align}
	- \frac{\log(1-z)}{z} & = \sum_{k=0}^\infty \frac{z^k}{k+1} \nonumber \\
	&= {\frac { \left( z-2 \right) }{ 2 \left( z-1
			\right) }} \sum_{n=0}^\infty \frac{1}{ \left( 2\,n+1 \right) {2n \choose n}}   \left( \frac{z^2}{z-1} \right)^n \label{z/bino(2n,n)} \\
	&= \frac{1}{3(z-1)} \, \sum_{n=0}^\infty \frac {\left( ((z+3) (2 z-3)) n + {z}^{2}+3\,z-6 \right) }{ \left( 3\,n+1 \right)  \left( 3
		\,n +2\right) {3n \choose n}} \left( \frac{z^3}{z-1} \right)^n \label{z/bino(3n,n)} \\
	&= \frac{1}{8(z-1)} \, \sum_{n=0}^\infty \frac {p(n) }{ \left( 4\,n+1 \right)  \left( 4
		\,n +2\right) (4n+3) {4n \choose n}} \left( \frac{z^4}{z-1} \right)^n, \label{z/bino(3n,n)}
\end{align}
where
{\small 
	\[
	p(n) = (9 z^3+12 z^2+16 z-64) n^2+(9 z^3+16 z^2+28 z-80) n+2 z^3+4 z^2+12 z-24.
	\]
}
\end{prop}

\noindent {\bf Remark.} By taking $z=-1,\frac{1}{2}, - \frac{1}{2}, \frac{1}{4}$, we recover \cite[Theorem 1.2]{Sun26}.

From the WZ seeds \eqref{WZ-abc'}--\eqref{WZS6} in Theorem~\ref{th-WZ}, we may derive more identities. We list some of them as follows.

Substituting  $a \to n, b \to  n, c \to  1/2+n, k \to 1+k+n$ into \eqref{WZ-abc'}, we derive that
\[
\sum_{n=1}^\infty {\frac {p(n)}{{n}^{3} \left( 3\,n-2 \right) ^{2}
	\left( -1+3\,n \right) ^{2} \left( -1+2\,n \right)  \ {5\,n\choose 3\,n}^{2} {2\,n\choose n} }}
=3 \pi^2,
\]
where
\begin{multline*}
p(n) = 574277\,{n}^{6}-1347180\,{n}^{5}+1274841\,{n}^{4}\\
-621294\,{n}^
	{3}+164024\,{n}^{2}-22184\,n+1200.
\end{multline*}

Taking $a = b= c = 1/2$ and $d =e= n$ in \eqref{WZ-abc11}, we derive that
\[
\sum_{n=0}^\infty (120n^2+34n+3) {\frac {{4\,n\choose 2\,n}  {2\,n\choose n} ^{4}}{{2}^{
		16\,n}}}
= \frac{2}{\pi} \sum_{n=0}^\infty (42n+5) \frac{{2n \choose n}^3}{2^{12n}}.
\]
By Guillera \cite[Identity 10]{Gu08}, the sum on the left-hand side equals $\frac{32}{\pi^2}$ and the sum on the right-hand side equals $\frac{16}{\pi}$ by Ramanujan \cite[Eq. (26)]{Ra14}.

Taking  $a = b= 1/2$,  $c=1/2+n$ and $d=e=n$ in \eqref{WZ-abc11}, we derive that
\begin{align*}
& \sum_{n=0}^\infty \frac{(20n^2-4n-1) n^2 (2n+1)}{2n-1} \frac{(-1)^n {2n \choose n}^5}{2^{12n}} \\
= & -\frac{1}{2 \pi} \sum_{n=0}^\infty (6n+1) \frac{{2n \choose n}^3}{2^{8n}} = - \frac{2}{\pi^2}.
\end{align*}
The last equality holds by Ramanujan \cite[Eq. (28)]{Ra14}, see also Baruah and Berndt \cite{BB10}.

Substituting  $a \to 1/2+n, b \to n, c \to 1, d \to 1+n, k \to 1+k+n$  into \eqref{WZS4}, we derive that
\[
\sum_{n=1}^\infty {\frac { \left( 224\,{n}^{4}+88\,{n}^{3}-66\,{n}^{2}-n+1 \right) 
}{ \left( 4\,n+1 \right)  \left( 2\,n+1
	\right) n \left( 2\,n-1 \right) ^{3} \left( 4\,n-1 \right) }} \frac{	\left( -1 \right) ^{n} {2n \choose n}}{{3n \choose n} 2^{4n}}
= 4 \log 2 - 3.
\]
Taking $a = 1/2-n, b = 1/2, c = 1/2, d = 1+n$  in \eqref{WZS4}, we derive that
\[
\sum_{n=1}^\infty {\frac {p(n)}{{n}^{5} \left( -1+3\,n \right) ^{2} \left( 3
	\,n-2 \right) ^{2}}}
{\frac {\left( -1 \right) ^{n} {2}^{16\,n} {3\,n\choose n} }{{6\,n
			\choose 3\,n}  {4\,n\choose 2\,n}^{3} {2\,n
			\choose n}^{2}}} = - 4032 \zeta(3),	
\]
where
\[
p(n) = 40752\,{n}^{6}-90240\,{n}^{5}+78696\,{n}^{4}-34564\,{n}^{3}+
8074\,{n}^{2}-954\,n+45.
\]

Substituting  $a,b,c,d \to \frac{1}{2}- n, e \to 1+n, k \to k+n$  into \eqref{WZS5}, we derive that
\[
\sum_{n=0}^\infty \frac{p(n) (-2^{24})^n }{(2n+1)^5 (6n+1)^5 (6n+5)^5 {6n \choose 3n}^5}
= \frac{189}{2}  \zeta(3),
\]
where
{
\small
\begin{multline*}
p(n) = 491287680\,{n}^{12}+3517115904\,{n}^{11}+11399133888\,{n}^{10} \\
+22100247360\,{n}^{9}+28522562400\,{n}^{8}+25791650640\,{n}^{7} \\
+16738310448\,{n}^{6}+7846312104\,{n}^{5}+2633313600\,{n}^{4} \\
+616193160\,{n}^{3}+95283770\,{n}^{2}+8729374\,n+357931.
\end{multline*}
}

Substituting $a,b,c,d,e \to \frac{1}{2}, f \to 1+n, k \to k+n$  into \eqref{WZS6}, we get
\begin{multline*}
F(n,k) = {\frac { \left( 8\,n+2\,k+1 \right)  \left( 3\,n+1+2\,k \right) }{
	\left( 4\,n+2\,k+1 \right) ^{5}}} \\
\times	\frac{\left( \frac{1}{2}-k \right)_n \left( \frac{1}{8}+\frac{k}{4} \right)_n \left( \frac{3}{8}+\frac{k}{4} \right)_n \left( \frac{5}{8}+\frac{k}{4} \right)_n \left( \frac{7}{8}+\frac{k}{4} \right)_n \left( \frac{1}{2}+k \right)_n^5 (1)_n^5}
{\left( \frac{1}{2} \right)_n^5 \left( \frac{1}{4}+\frac{k}{2} \right)_n^5 \left( \frac{3}{4}+\frac{k}{2} \right)_n^5}.
\end{multline*}
We thus derive that
\[
g(n) = \lim_{k \to \infty} G(n,k) = -\frac{205 n^2+250 n+77}{512} \frac{ (1)_n^5}{ \left( \frac{3}{2} \right)_n^5} \left( -\frac{1}{1024} \right)^n,
\]
and hence
\begin{align*}
\sum_{n=0}^\infty G(n,0)  & =
\frac{1}{512} \sum_{n=0}^\infty \frac{p(n) (8n+1) }{(4n+1)^5 (4n+3)^5} \frac{ {8n \choose 4n}}{{4n \choose 2n}^4 {2n \choose n}^2} \\
& = \sum_{k=0}^\infty F(0,k) + \sum_{n=0}^\infty g(n) \\
& = \sum_{k=0}^\infty \frac{1}{(2k+1)^3} - \frac{1}{512} \sum_{n=0}^\infty (205n^2+250n+77)  \frac{ (1)_n^5}{ \left( \frac{3}{2} \right)_n^5} \left( -\frac{1}{1024} \right)^n \\
& = \frac{7}{8} \zeta(3) - \frac{1}{8} \zeta(3)  = \frac{3}{4} \zeta(3),
\end{align*}
where
\begin{multline*}
	p(n) = 1397760\,{n}^{6}+5727104\,{n}^{5}+9673968\,{n}^{4}+8600034\,{n}^{3}\\
	+
	4227785\,{n}^{2}+1083660\,n+112149
\end{multline*}
and we use the following identity due to Amdeberhan and Zeilberger \cite{AZ97}:
\[
\sum_{n=0}^{\infty} (-1)^n \frac{n!^{10} (205n^2 +250n+77)}{64 (2n+1)!^5} = \zeta(3).
\]

\section{$q$-analogues}
Recall that a function $F(n_1, \ldots,n_r)$ is  {\it $q$-hypergeometric} in the variables $n_1,\ldots,n_r$ if all the ratios 
\[
\frac{F(n_1+1,n_2,\ldots,n_r)}{F(n_1,n_2,\ldots,n_r)},\ \frac{F(n_1,n_2+1,n_3,\ldots,n_r)}{F(n_1,n_2,\ldots,n_r)}, \ \ldots, \frac{F(n_1,\ldots,n_{r-1}, n_r+1)}{F(n_1,n_2,\ldots,n_r)}
\]
are rational functions of $q^{n_1}, q^{n_2},\ldots, q^{n_r}$.  A pair of $q$-hypergeometric terms $(F(n,k), G(n,k))$ forms a {\it $q$-WZ pair} if \eqref{WZ} holds. In this case, we call $G(n,k)$ the {\it $q$-WZ mate} of $F(n,k)$. A $q$-hypergeometric term $F_0(k,a,b,\ldots)$ is a {\it $q$-WZ seed} if for all integers $K,A,B,\ldots$, 
\[
F(n,k) = F_0(Kn+k_0+k, An+a,Bn+b,\ldots)
\]
has a $q$-WZ mate. 

Suppose that when we apply the $q$-Zeilberger algorithm to a $q$-hypergeometric term $F(k,a,b,\ldots)$ with shifts in $a$ and obtain the telescoped recurrence of order one
\[
p_0(a,b,\ldots) F(k, a,b,\ldots) + p_1(a,b,\ldots) F(k, a+1,b,\ldots) = \Delta_k G(k,a,b,\ldots),
\]
where $p_0$ and $p_1$ are polynomials in $q^a$ and independent of $q^k$. Factoring $p_0$ and $p_r$, we can write
\begin{equation*}
p_0 = c_0 \cdot (q^a)^s \cdot  (1 - u_1 q^a)  \ldots (1 - u_l q^a), \quad p_1 = c_1 \cdot (q^a)^t \cdot (1 - v_1 q^a)  \ldots (1 - v_m q^a).
\end{equation*}
By setting
\begin{equation*}
F_a(k,a,b,\ldots) =   \left(-\frac{c_1}{c_0} \right)^a \left( q^{a \choose 2} \right)^{t-s} \frac{\prod_{i=1}^l (u_i q^a;q)_\infty}{\prod_{j=1}^m (v_j q^a; q)_\infty} \cdot F(k, a, b,\ldots),
\end{equation*}
where $(a;q)_\infty = \prod\limits_{i=0}^\infty (1-a q^i)$, 
we find that $F_a(k,a,b,\ldots)$ has a $q$-WZ mate.  We remark that one can make the replacement \begin{equation}\label{qrep}
q^{a \choose 2} \to \frac{(-1)^a}{(q^a;q)_\infty (q^{1-a};q)_\infty},
\end{equation}
which is nicer when $a$ is not an integer. 

By the above construction, we obtain  $q$-analogues of Theorem~\ref{th-WZ}.
Here and in what follows, we use the shorten notation
\[
\left[
\begin{array}{c}
a_1,\ldots, a_r \\
b_1,\ldots, b_s
\end{array}; q
\right]_\infty = \frac{\prod_{i=1}^r (a_i;q)_\infty}{\prod_{j=1}^r (b_j;q)_\infty}.
\]
\begin{theorem}\label{qWZ-seed}
The following functions are $q$-WZ seeds.
\begin{itemize}
	
	\item[{\rm (1')}]
	\begin{equation}\label{WZ-abz-q}
		z^{b+k} 	\left[
		\begin{array}{c}
			q^{b+k} \\ q^{a+k} 
		\end{array}; q
		\right]_\infty \cdot 
		\left[
		\begin{array}{c}
		 q^{a-b+1} \\  z q^{a-b+1}
		\end{array}; q
		\right]_\infty,
	\end{equation}
	
	\item[{\rm (2')}]
	\begin{multline}\label{th-abcd-q}
		(-1)^{c+d} q^{k+{c+1 \choose 2}+{d+1 \choose 2}-cd} \left[
		\begin{array}{c}
			q^{c+k},\ q^{d+k} \\
			q^{a+k},\ q^{b+k}
		\end{array}; q
		\right]_\infty \\
		\times 
		\left[
		\begin{array}{c}
			q^{a-c+1},\ q^{a-d+1},\ q^{b-c+1},\ q^{b-d+1} \\
			q^{a+b-c-d+2}
		\end{array}; q
		\right]_\infty
	\end{multline}
	
	\item[{\rm (3')}]
	\begin{multline}\label{WZ-abc'-q}
		(-1)^k {q}^{2c k+ {2c \choose 2} +4\,bc}
		\left[
		\begin{array}{c}
			{q}^{2\,c+2\,a+k},\ {q}^{2\,b+2\,c+k} \\
			{q}^{2\,a+k},\ {q}^{2\,b+k}
		\end{array}; q
		\right]_\infty  \\[10pt]
		\times \left[
		\begin{array}{c}
			-{q}^{2\,c} \\
			{q}^{2\,c}
		\end{array}; q
		\right]_\infty
		\cdot \left[
		\begin{array}{c}
			{q}^{2\,a-2\,b-2\,c+1} \\ {q}^{2\,a-2\,b+2\,c+1}
		\end{array}; q^2
		\right]_\infty,
	\end{multline}
	
	\item[{\rm (4')}]
	\begin{multline}\label{WZ-abc11-q}
		\left( -1 \right) ^{d+e}{q}^{k + {d+1 \choose 2} + {e+1 \choose 2} - de}
		\left[
		\begin{array}{c}
			{q}^{d+k},\ {q}^{e+k},\ {q}^{a+b+c-d-e+2+k} \\
			{q}^{a+k},\ {q}^{b+k},\ {q}^{c+k}
		\end{array}; q
		\right]_\infty \\[10pt]
		\times
		\left[
		\begin{array}{c}
			{q}^{a-d+1},\ {q}^{a-e+1},\ {q}^{b-d+1},\ {q}^{b-e+1},\
			{q}^{c-d+1},\ {q}^{c-e+1} \\
			{q}^{a+b-d-e+2},{q}^{a+c-d-e+2},{q}^{b+c-d-e+2}
		\end{array}; q
		\right]_\infty,
	\end{multline}
	
	\item[{\rm (5')}]
	\begin{multline}\label{WZS4-q}
		(-1)^{a+b+c} {q}^{ \left( 1+a-b-c \right) k+ {a+1 \choose 2} +{2a-b-c+1 \choose 2}  - bc } \\[5pt]
		\times	\left( 1+{q}^{a+k} \right) \cdot  \left[
		\begin{array}{c}
			{q}^{1+k},\ {q}^{2a-b+1+k},\ {q}^{2a-c+1+k} \\
			{q}^{2\,a+k},\ {q}^{b+k},\ {q}^{c+k}
		\end{array}; q
		\right]_\infty \\[10pt]
		\times
		\left[
		\begin{array}{c}
			q^a,\ q^b,\ q^c,\ {q}^{b-a},\ {q}^{c-a},\ {q}^{b+c-2a} \\
			{q}^{b+c-a}
		\end{array}; q
		\right]_\infty,
	\end{multline}
	
	\item[{\rm (6')}]
	\begin{multline}\label{WZS5-q}
		\left( -1 \right) ^{a+b+c+d+e} q^{(2e-a-b-c-d+1)k  + {a-e \choose 2}+ {b-e \choose 2}+ {c-e \choose 2}+ {d-e \choose 2} - {e+1 \choose 2}} \\[5pt] 
		\times (1-{q}^{2\,k+e}) \cdot 
		\left[
		\begin{array}{c}
			{q}^{e-a + 1 +k},\  {q}^{e-b + 1 +k},\  {q}^{e-c + 1 +k},\ {q}^{e-d + 1 +k} \\
			{q}^{a+k},\ {q}^{b+k},\ {q}^{c+k},\ {q}^{d+k}
		\end{array}; q
		\right]_\infty \\[5pt] 
		\times
		\left[
		\begin{array}{c}
			{q}^{a+b-e},\ {q}^{a+c-e},\ {q}^{a+d-e},\ {q}^{b+c-e},\ {q}^{b+d-e},\ {q}^{c+d-e} \\
			q^{a+b+c+d-2e}
		\end{array}; q
		\right]_\infty,	
	\end{multline}
	
	\item[{\rm (7')}]
	{\small 
		\begin{multline}\label{WZS6-q}
			\left( -1 \right) ^{f}{q}^{k+ {a+1 \choose 2} + {b+1 \choose 2} + {c+1 \choose 2} + {d+1 \choose 2} + {e+1 \choose 2} - {f+1 \choose 2} - {2f-a-b-c-d-e \choose 2}}  \\[5pt] 
			\times (1-{q}^{f+2\,k}) \cdot \left[
			\begin{array}{c}
				{q}^{f-a+1+k},\ {q}^{f-b+1+k},\ {q}^{f-c+1+k},\ {q}^{f-d+1+k},\ {q}^{f-e+1+k}
				\\
				{q}^{a+k},\ {q}^{b+k},\ {q}^{c+k},\ {q}^{d+k},\ {q}^{e+k}
			\end{array}; q
			\right]_\infty \\[5pt] 
			\times 
			\left[
			\begin{array}{c}
				{q}^{a+b+c+d+e-2f+k}
				\\
				q^{3f-a-b-c-d-e+1+k}
			\end{array}; q
			\right]_\infty \cdot
			\left[
			\begin{array}{c}
				{q}^{a+b-f},\ {q}^{a+c-f},\ {q}^{a+d-f},\ {q}^{a+e-f}
				\\
				{q}^{a+b+c+d-2f},\ {q}^{a+b+c+e-2f}
			\end{array}; q
			\right]_\infty \\[5pt] 
			\times
			\left[
			\begin{array}{c}
				{q}^{b+c-f},\ {q}^{b+d-f},\ {q}^{b+e-f},\ {q}^{c+d-f},\ {q}^{c+e-f},\ {q}^{d+e-f}
				\\
				{q}^{a+b+d+e-2f},\ {q}^{a+c+d+e-2f},\ {q}^{b+c+d+e-2f}
			\end{array}; q
			\right]_\infty.
		\end{multline}
	}
\end{itemize}
\end{theorem}

\noindent {\bf Remark.} When $q \to 1$, \eqref{WZS4-q} turns into the special case of \eqref{WZS4} with $a \to 2a$ and $d \to a+\frac{1}{2}$. The remaining $q$-WZ seeds become their corresponding WZ seeds.

Notice that we may make the replacement 
\[
(q^a;q)_\infty \to  \frac{(-1)^a q^{-{a \choose 2}}}{(q^{1-a};q)_\infty}
\]
in the expressions of $q$-WZ seeds.  Utilizing this replacement, we recover $q$-WZ seeds given by Au \cite{Auc}:
\begin{itemize}
\item 
By setting $e=a+b+c+1-d$, \eqref{WZ-abc11-q} is equivalent to {\tt Balanced3F2$_q$}.
\item 
\eqref{WZS4-q} is equivalent to  {\tt Dixon3F2$_q$}.
\item 
By setting $e=a$, \eqref{WZS5-q} is equivalent to {\tt Dougall5F4$_q$}.
\item 
By setting $f=a$, \eqref{WZS6-q} is equivalent to {\tt Dougall7F6$_q$}.
\end{itemize}

Noting that
\[
(q^a;q)_k = (-q^a)^k q^{k \choose 2} (q^{-a};q^{-1})_k.
\]
We may replace $q$ with $q^{-1}$ in a hypergeometric term to derive a new WZ seed. For instance, when we replace $q$ with $q^{-1}$ and $a,b,c,d$ with $-a,-b,-c,-d$ respectively in 
\[
q^k \left[
\begin{array}{c}
q^{a},\ q^b
\\
q^c,\ q^d
\end{array}; q
\right]_k,
\]
we will obtain
\[
q^{(c+d-a-b-1) k} \left[
\begin{array}{c}
q^{a},\ q^b
\\
q^c,\ q^d
\end{array}; q
\right]_k.
\]
This leads to the $q$-WZ seed
\begin{multline}
(-1)^{c+d}  q^{(c+d-a-b-1)k + {c \choose 2} + {d \choose 2} - ab} \\
\times 
\left[
\begin{array}{c}
	q^{c+k},\ q^{d+k} \\
	q^{a+k},\ q^{b+k}
\end{array}; q
\right]_\infty
\left[
\begin{array}{c}
	q^{a-c+1},\ q^{a-d+1},\ q^{b-c+1},\ q^{b-d+1} \\
	q^{a+b-c-d+2}
\end{array}; q
\right]_\infty, \label{q-abcd}
\end{multline}
which differs from \eqref{th-abcd-q} only by a factor of power of $q$. When $d=1$, the resulting $q$-WZ seed is equivalent to {\tt Gauss2F1$_q$}  given by Au \cite{Auc}.

We conclude this section with a simple application. 
We use the above $q$-WZ seed to give a $q$-analogue of a formula due to Guillera \cite{Gu08}
\[
\sum_{n=0}^\infty \frac{(3n+2) 2^{4n}}{(2n+1)^3 {2n \choose n}^3} = \frac{\pi^2}{4}.
\]
By taking $a=b=\frac{1}{2}, c=d=\frac{3}{2}+n$ in \eqref{q-abcd}, we get
\[
F(n,k) = q^{(2n+1) k+n+{n}^{2}}  \frac{(q^{\frac{3}{2}+n+k};q)_\infty^2 (q^{2n+1};q)_\infty}
{(q^{\frac{1}{2}+k};q)_\infty^2 (q^{n+1};q)_\infty^4},
\]
and
\[
G(n,k) = \frac {2{q}^{2n+\frac{3}{2}+k} - q^{n+1} - 1}{
\left( {q}^{2n+1} -1 \right)  \left( {q}^{n+1}+1
\right) } F(n,k).
\]
This leads to 
\[
\sum_{k=0}^\infty \frac{q^k}{(1-q^{k+\frac{1}{2}})^2} = \sum_{n=0}^\infty
\frac{(1-2q^{2n+\frac{3}{2}} + q^{n+1}) q^{n^2+n}}{(1-q^{2n+1})(1+q^{n+1})}
\frac{(q;q)_n^4}{(q^{\frac{1}{2}};q)_{n+1}^2 (q;q)_{2n}},
\]
which is the same as the one given in \cite{HKS}.

\end{document}